\documentclass[10pt,preprint]{imsart}

\RequirePackage[OT1]{fontenc}
\RequirePackage{amsthm,amsmath}
\RequirePackage[colorlinks,citecolor=blue,urlcolor=blue]{hyperref}

\usepackage[square,compress,comma, numbers]{natbib}
\usepackage{amsfonts,mathtools}

\def\nu{\theta}
\def\tX{\theta_X}

\startlocaldefs
\numberwithin{equation}{section}
\theoremstyle{plain}

\endlocaldefs

\def\IA{\mathcal{I}_{fm}}
\def\IAfE{\mathcal{I}_{fe}}

\allowdisplaybreaks[4]

\def\tX{\widetilde{X}}

\usepackage{amssymb}
\usepackage{color}
\newcommand{\abs}[1]{\left\lvert #1 \right\rvert}

\DeclarePairedDelimiterXPP\pk[1]{\mathbb{P}}\{ \}{}{ #1}
\DeclarePairedDelimiterXPP\E[1]{\mathbb{E}}\{ \}{}{	#1}

\def\FRE{\mbox{Fr\'{e}chet }}

\usepackage{xparse}

\NewDocumentCommand{\ceil}{s O{} m}{%
	\IfBooleanTF{#1} 
	{\left\lceil#3\right\rceil} 
	{#2\lceil#3#2\rceil} 
}
\NewDocumentCommand{\floor}{s O{} m}{%
	\IfBooleanTF{#1} 
	{\left\lfloor#3\right\rfloor}
	{#2\lfloor#3#2\rfloor}
}

\newcommand{\norm}[1]{\lVert #1 \rVert}

\definecolor{c20}{rgb}{0.,0.7,0.}
\definecolor{c30}{rgb}{0.,0.,1.}
\definecolor{c40}{rgb}{1,0.1,0.7}
\definecolor{c50}{rgb}{1,0,0}
\definecolor{c60}{rgb}{1,0.9,0.1}
\definecolor{c70}{rgb}{0.50,1.00,0.00}

\def\cL#1{{\textcolor{c30}{#1}}}
\def\cL#1{#1}

\def\EHF#1{{\textcolor{c50}{#1}}}
\def\EHF#1{#1}

\numberwithin{equation}{section}
\newtheorem{theo}{Theorem}[section]
\newtheorem{sat}[theo]{Proposition}
\newtheorem{de}[theo]{Definition}
\newtheorem{lem}[theo]{Lemma}

\newtheorem{example}[theo]{Example}
\newtheorem{korr}[theo]{Corollary}
\newtheorem{remark}[theo]{Remark}

\numberwithin{equation}{section}

\newcommand{\prooftheo}[1]{ \textsc{Proof of Theorem} \ref{#1}:}

\newcommand{\prooflem}[1]{\textsc{Proof of Lemma} \ref{#1}:}

\newcommand{\QED}{\hfill $\Box$}

\newcommand{\COM}[1]{}

\def\IF{\infty}

\newcommand{\R}{\mathbb{R}}
\newcommand{\inr}{\in \R}

\newcommand{\BQN}{\begin{eqnarray}}
\newcommand{\EQN}{\end{eqnarray}}
\newcommand{\BQNY}{\begin{eqnarray*}}
	\newcommand{\EQNY}{\end{eqnarray*}}
\def\ldot{, \ldots,}

\newcommand{\limit}[1]{\lim_{#1 \to   \infty}}

\def\todis{\overset{d}\rightarrow}

\def\CE#1{\textcolor{c30}{#1}}

\def\bqny#1{\begin{eqnarray*} #1 \end{eqnarray*}}
\def\bqn#1{\begin{eqnarray} #1 \end{eqnarray}}

\newcommand{\BS}{\begin{sat}}
	\newcommand{\ES}{\end{sat}}
\newcommand{\BT}{\begin{theo}}
	\newcommand{\ET}{\end{theo}}
\newcommand{\BK}{\begin{korr}}
	\newcommand{\EK}{\end{korr}}

\newcommand{\BEX}{\begin{example}}
	\newcommand{\EEX}{\end{example}}

\newcommand{\BD}{\begin{de}}
	\newcommand{\ED}{\end{de}}

\newcommand{\BIT}{\begin{itemize}}
	\newcommand{\EIT}{\end{itemize}}
\newcommand{\BDI}{\begin{description}}
	\newcommand{\EDI}{\end{description}}

\newcommand{\BRM}{\begin{remark}}
	\newcommand{\ERM}{\end{remark}}

\newcommand{\BEL}{\begin{lem}}
	\newcommand{\EEL}{\end{lem}}

\def\SSZ{ \mathcal{S}(Z)}
\def\SST{ \mathcal{S}(\Theta)}
\def\JJZ{ \mathcal{I}(Z)}
\def\JJT{ \mathcal{I}(\Theta)}
\def\JJY{ \mathcal{I}(Y)}

\newcommand{\nelem}[1]{{Lemma \ref{#1}}}

\newcommand{\netheo}[1]{{Theorem \ref{#1}}}

\newcommand{\equaldis}{\stackrel{{d}}{=}}

\def\Z{\mathbb{Z}}
\def\inn{\in \mathbb{N}}

\newcommand\ind[1]{\mathbb{I}{( #1)}}

\def\TT{\Z^d}

\begin{document}

\begin{frontmatter}
\title{On  Extremal Index of Max-Stable Random Fields}

\begin{aug}
\author{\snm{Enkelejd Hashorva }\ead[label=e1]{enkelejd.hashorva@unil.ch}}


\affiliation{University of Lausanne}

\address{Department of Actuarial Science, University of Lausanne\\
	Chamberonne 1015 Lausanne, Switzerland\\
\printead{e1}\\
\phantom{E-mail:\ }}

\end{aug}
\def\tX{\theta_X}

\begin{abstract}
	For a given stationary max-stable random field $X(t),t\in \TT$ the corresponding generalised Pickands constant coincides with the classical extremal index $\tX \in [0,1]$.  
	In this contribution we discuss  necessary and sufficient conditions for $\tX$ to be 0,  positive 
	or equal to 1  and  also show  that $\tX$  is equal to the so-called block extremal index. 
	Further, we consider some general functional  indices  of $X$  and prove  that for a large class of  functionals they coincide with $\tX$. 
 \end{abstract}

\begin{keyword}[class=MSC]
\kwd[Primary ]{60G15}
\kwd[; secondary ]{60G70}
\end{keyword}

\begin{keyword}
\kwd{Max-stable random fields}
\kwd{Brown-Resnick random fields}
\kwd{Pickands constants}
\kwd{classical extremal index}
\kwd{block extremal index}
\kwd{functional  index}

\end{keyword}

\end{frontmatter}

\section{Introduction}
The connection between Pickands constant and extremal index of stationary max-stable  Brown-Resnick random fields (rf's) has been initially pointed out in \citep{DM}. Calculation of Pickands constants for a general  stationary max-stable rf $X(t),t\in \TT $ has been later dealt with in  \cite{Htilt}. Previous investigations  concerned  with the calculation of extremal index in the context of max-stable processes are \citep{stoev2010max,MR2453345,MIKYuw,KDEH1}. 
Recent research in \cite{BP,WS,SojaA,SojaB} has shown, contrary to the prevailing intuitions, that there are certain  subtilities (if $d>1$) when dealing with stationary multivariate regularly varying rf's (see e.g., \cite{MR3764610} for the definition)  and the calculation of their extremal indices. 
Influenced by  the findings of \cite{Davis},  several formulas for  extremal indices of stationary  regularly varying time series have appeared in the literature, see e.g., \cite{Hrovje} and the references therein. Various (less well-known) formulas have been discovered  also for Pickands constants in contributions unrelated to time series modelling. For instance in sequential analysis and statistical applications \cite{Sigm4, Sigm5} and extremes of  random fields \cite{BYak, Pit19} just to mention a few. 
For large classes of Gaussian rf's 
extremal indices   have been discussed  in \cite{DavisF,MR3336846,MR3813912}, see also \cite{MR2225054,MR3226001} for non-Gaussian cases and related results. \\ 

Without loss of generality,   we shall focus on the class of max-stable rf's with $\FRE$ marginals. Since these  are limiting rf's,  see e.g., \cite{klem},  our formulas for their extremal indices are valid (with obvious modifications) also for the candidate extremal index of more general stationary  regularly varying rf's (see  \cite{Hrovje} for recent findings). Studying max-stable rf's, instead of these more general rf's is also justified by   \nelem{prop2} stated in Section 2 \EHF{and Remark \ref{remF} $iii)$}.\\
In view of the well-known de Haan characterisation given in \cite{deHaan}, the rf $X$ with non-degenerated marginal distributions corresponds to some non-negative spectral rf $Z(t),t\in \TT$ having  the following representation (in distribution) 
\bqn{\label{eq1}
	X(t) =  \max_{i\ge 1} \Gamma_i^{-1/\alpha} Z_i(t), \quad t\in \TT,
}
where $\Gamma_i= \sum_{k=1}^i Q_k$ with $Q_k, k\ge 1$ unit exponential random variables \EHF{(rv's)} independent of $Z_i$'s which are independent copies of $Z$.\\
\EHF{Clearly, $Z$ is not unique since also $\tilde Z(t)=R Z(t),t\in \TT$ is a spectral rf for $X$, provided that $R$ is a non-negative rv independent of $Z$ such that $\E{R^\alpha}=1$. Note that if for some $h\in \TT$ we have $Z(h)=1$ almost surely, then in view of Balkema's lemma (stated in \cite{HaanPickands}[Lem 4.1]) any spectral rf $\tilde  Z$ of $X$ has the same law as $Z$. }\\
 We shall assume without loss of generality that for some $\alpha \in (0,\IF)$ 
\bqn{ \label{elit}
	\pk*{\max_{t\in \TT} Z(t)>0}=1, \quad 
	\E*{ Z^\alpha (t)}=1, \quad t\in \TT. 
}
\nelem{lemX} in Appendix  shows how to construct a spectral rf $Z$ such that the first assumption in \eqref{elit}  holds. Note that $\E*{Z^\alpha(t)} =1$ 
 implies that $X(t)$  has   $\alpha$-\FRE distribution function $e^{-x^{-\alpha}},x>0$. This is no restriction since we are interested in stationary max-stable rf's.  
As in \cite{Htilt} define the  Pickands constant (when the limit exists) with respect to the spectral rf $Z$  by 
\def\HWD{\mathcal{H}}
\bqn{\label{aha}
	\HWD = \limit{n} \frac 1 {n^d} \E*{ \max_{t\in [0,n]^d\cap \Z^d } Z^{\EHF{\alpha}}(t) }\le 
	\limit{n} \frac 1 {n^d} \sum_{t\in [0,n]^d\cap \Z^d }  \E*{ Z^\alpha(t) }\le 1 .
} 
Since the finite dimensional distributions (fidi's) of $X$ can be calculated explicitly (see \eqref{fidis} below), if $\HWD$ exists, then   
\bqn{\label{tiz}
	\quad \quad 
	\pk*{ \max_{t\in [0,n]^d \cap \Z^d } X(t) \le n^d x}= e^{-  \frac{1 }{n^d}\E*{ \max_{t\in [0,n] ^d  \cap \Z^d} Z^\alpha(t)/x^\alpha  } }\to 
	e^{-   \HWD/ x^\alpha}
}
as $n\to \IF $ is valid for   all $x>0$.
\def\tX{\theta_X}

\noindent
As argued  in \cite{DM} and \cite{Htilt,debicki2017approximation}   the sub-additivity of maximum functional implies that    $\HWD$ is  well-defined and finite, provided that $X$ is stationary. Consequently, in view of \eqref{tiz} the extremal index  (or using the terminology of \cite{WS},  the classical extremal index) of the stationary max-stable rf $X$ (denoted below by $\tX$)  \EHF{always exists, does not depend on the particular spectral rf $Z$ but on the law of the rf $X$  and is given by} 
\bqn{\label{AVI}
	\tX= \HWD \in [0,1]. 
}	
\EHF{In the special case  
	\bqn{ \label{lfever} X(t) =V_t, \quad t\in \TT,
	}
 where $V_t$'s are independent $\alpha$-\FRE rv's we have  $\tX=1$. We shall show that this is the only max-stable rf with unit \FRE marginals satisfying  $\tX=1$. Using this fact and  \nelem{prop2}  we can construct a  spectral rf $Z$ for $X$, see Remark \ref{nms} $iii)$.}

 Hereafter  
we shall assume  for simplicity that the max-stable rf $X$ has unit \FRE marginal distributions,  i.e., below we shall consider the case   
$$ \alpha=1.$$

If the spectral rf $Z$ is not easy to determine or  $X(t), t\in \TT$ is  stationary  but not max-stable,  commonly the  block extremal index  (denoted below by $\widetilde{\tX}$)  is utilised in various applications related to extreme value analysis. Assuming for simplicity that $X$ has unit \FRE marginals, it is  defined by (see \cite{MR2444320, WS})
\bqn{ \widetilde{\tX}:=  \limit{n}  \frac{\pk{ \max_{ 0 \le   i \le r_n, i\in \TT } X(i)> n  \tau} }{ \prod_{j=1}^d r_{nj}  
		\pk{ X(0)> n \tau}}
	\label{bloc}
}
for any $\tau>0$ and any  sequence $r_n\in \TT,n \ge1$ with non-decreasing  integer-valued components $r_{nj}, j\le d$ such that $\limit{n} r_{nj}= \limit{n} n/r_{nj}^d =\IF$ for any $j\le d$. In our setting we do not need to put the last restriction. 
In \eqref{bloc}  $i \le r_n$ is interpreted component-wise, i.e., $i_j \le r_{n_j}$ for all $j\le d$ components of $i$ and $r_n$, respectively.

Next, we define    functional  indices $\nu_{X,F}$ of $X$ by 
$$ \nu_{X,F}=\E*{ Z(0) F(Z)} \EHF{\in [0,1]},$$
where $F: E \mapsto [0,1]$ is a measurable functional with respect to the product \cL{$\sigma$-field}   $\mathcal{E}$  on $E:=[0,\IF)^{\TT}.$ \\
\EHF{As mentioned above different choices of $Z$ for $X$ are possible.} In order to make the definition of $\nu_{X,F}$ independent of the choice of   $Z$ and thus only dependent on the law of $X$, we shall also require that $F$ is $0$-homogeneous, i.e., $F(cf) =F(f)$ for any $c>0, f\in E$. Indeed, under this assumption 
we have that 
$$
\nu_{X,F}= \E*{ Z(0) F(Z/Z(0))}= \E*{ F(\Theta_{0})}  ,
$$   
where  the rf $\Theta_{h}$  is  defined by (hereafter $\ind{\cdot}$ denotes the indicator function)
\bqn{\label{eqtheta} \pk{\Theta_{h} \in A} &=&  \E*{ Z(h) \mathbb{I}(Z/Z(h)\in A )	}, \quad \forall A \in   \mathcal{E}.
}   
\EHF{It is known that for any $h\in \TT$ the law of $\Theta_h$  does not depend on the particular choice of the spectral rf $Z$ and can be directly determined by $X$. 
In the case that for a spectral rf $Z$ of $X$ we have that $Z(h)>0$ almost surely, this fact follows from Balkema's lemma. The proof for the general case follows from \cite{Htilt}[Lem A.1], or from \cite{StoevWangSPL}[Thm 1.1] and  \cite{MolchanovBE}[Thm 2].  
Consequently,	the functional  index $ \nu_{X,F}$ depends only on the law of $X$. Note that for the definition of $\nu_{X,F}$ no stationarity of $X$ is assumed. }\\
 
\EHF{
	It is well-known that a max-stable rf $X$ with \FRE marginals is  a multivariate regularly varying rf. For general multivariate regularly varying rf's which are not max-stable,  there is no spectral process $Z$ as in our case of max-stable $X$ and therefore the  rf's 
	$\Theta_h, h \in \TT$ are defined via  a conditional limit, see e.g., \cite{klem,SegersEx17} and \eqref{convT} below. The key advantage in the framework  of max-stable rf's is that $\Theta_h$ is  directly obtained by tilting a given  spectral rf $Z$.
 }

 At this point two natural questions for a given stationary max-stable rf $X$ arise: \\
{\bf Question 1:}  What is the relation between $\tX$ and $\widetilde{\tX}$? \\ 
 
 \noindent
 {\bf Question 2:}  For what $F$ is the functional  index $\nu_{X,F}$  equal to $\tX$?\\
 
 \noindent

In this contribution we show that we simply have $\tX= \widetilde{\tX}$ and then describe a large class of  functionals   $F$ such that $\tX= \theta_{F,X}$. Further, we consider in some detail the cases $\tX=0$ and $\tX=1$. 
 
Brief organisation of the rest of the paper:  In the next section we discuss some basic properties of the  rf's $\Theta_h, h\in \TT$
 and then show how to construct a \EHF{stationary} max-stable rf $X$ from a given   rf $\Theta^*$ which in turn is necessary equal in law with $\Theta_0$. In Section 3 we claim that $\tX=\widetilde{\tX}$ for any  stationary max-stable rf's $X$.   Additionally, we give equivalent conditions that guarantee $\tX>0$ or $\tX=0$ and then present several formulas for $\tX$. Section 4 is concerned with the anti-clustering condition whereas Section 5 displays  some  examples. All the proofs are relegated to Section 6   which is followed by an Appendix.

\section{Preliminaries}
Unless otherwise specified we shall consider  below a max-stable rf $X(t),t\in \TT$  as in the Introduction with spectral rf $Z$ such that  $\E*{Z(t)}=1, t\in \TT$. Hence $X(t)$ has unit \FRE distribution $e^{-1/x},x>0$. We shall discuss first  the case that  $X$ is   non-stationary.

\subsection{General max-stable $X$}
The importance of the rf's $\Theta_h,h \in \TT$ defined in \eqref{eqtheta} relates to the following conditional convergence results. Namely,  
in view of \cite{Htilt}[Lem 2.1, A.1 \&  Rem 6.4] or by \cite{klem}[Lem 3.5]  we have that the convergence in distribution 
\bqn{\label{convT}
	X (t)/X(h) \Bigl \lvert (X(h) > u ) \todis \Theta_h(t), \quad t\in \TT, \\
	\label{convY}
	u^{-1} X (t) \Bigl \lvert (X(h) > u ) \todis Y_h(t), \quad t\in \TT
}
hold as $u\to \IF$ in the product topology of $E=[0,\IF)^{\Z^d}$, where $\Theta_h$ is defined in \eqref{eqtheta} 
and 
	$$Y_h(t)= R \Theta_h(t),\quad t\in \TT,$$
	with $R$  an $\alpha$-Pareto rv with survival function $x^{-\alpha} , x\ge 1$ independent of any other random element (recall that we consider $\alpha=1$ for simplicity).\\
 \EHF{If for a given  max-stable rf $X$ if a spectral rf $Z$ is known, it is often simpler  to determine the law of  $\Theta_h$ directly via \eqref{eqtheta} than deriving it from \eqref{convT}. In particular}, if \EHF{$\pk{Z(h)=1}=1$},  then the following equality in law   
\bqn{ \label{fle}
	\Theta_h \equaldis  Z
}	 
is valid. Below we determine the fidi's of  $Y_h $ in terms of $Z$ and  $\Theta_h $. 
\BEL \label{theo1} 
 For any $ h,t_i \in \TT,  x_i \in (0,\IF), i\le n$ we have 
 \bqn{ 
	\pk{ Y_{\EHF{h}}(t_1 ) \le x_1  \ldot Y_{\EHF{h}}(t_n) \le x_n} =
	\E*{\max\Bigl(1, \max_{1 \le i \le n } \frac{\Theta_{\EHF{h}} (t_i)}{x_i} \Bigr)- 	\max_{1 \le i \le n } \frac{\Theta_{\EHF{h}}(t_i)}{x_i}} \label{eB} \\   
=	\E*{ \max\Bigl( Z(h),  \max_{1 \le i \le n } \frac{Z (t_i)}{x_i} \Bigr)- 	\max_{1 \le i \le n } \frac{Z(t_i)}{x_i}}. \notag
} 
\EEL

\BRM   
For the case of the stationary Brown-Resnick model \eqref{eB} is stated in \cite{WS}[Prop 6.1] for $h=0$.\\
 \ERM

\subsection{Stationary max-stable $X$}
In view of  \cite{Htilt}[Thm 6.9] the max-stable rf $X(t),t\in \TT$ with unit \FRE marginals is stationary, if and only if 
\bqn{\label{eqDo}
	\E*{ Z (h) F(Z)} = \E*{Z (0) F(B^h Z)}, \quad \forall h\in \TT}
is valid for any  measurable function $F:  E \mapsto [0,\IF]$ which is $0$-homogeneous. 
Here $B$ is the shift-operator so that $B^h Z(\cdot) = Z(\cdot - \EHF{h}),h\in \TT$. Note  that for the stationary Brown-Resnick model the claim in \eqref{eqDo} is first formulated in \cite{DM}[Lem 5.2].  \\
\EHF{For notational simplicity we shall omit the subscript 0 and write simply $\Theta$ and $Y$ instead of $\Theta_0$ and $Y_0$, respectively;
	in our notation  the origin of $\R^k,k\inn$  is denoted by 0.}\\
  \EHF{ In view of \cite{Htilt}[Thm 4.3] condition  \eqref{eqDo} is equivalent with the following equality in law 
	$$\Theta_h \equaldis B^h \Theta$$
	 valid  for any $h\in \TT$.}\\
 Yet another equivalent formulation of condition  \eqref{eqDo}  stated for the rf $\Theta$ is 
\bqn{\label{eqDo2}
	\E*{{  \Theta (h)} F({ \Theta})}  
	&=& \E*{  F( B^{h} {\Theta }) \mathbb{I}( B^h \Theta(0) \not=0)}, \quad \forall h\in \TT
}
valid again for all measurable functionals $F$ as above, see e.g., \cite{BP,klem}.\\
\EHF{We note in passing that  with the same arguments as in  \cite{klem} it can be shown that \eqref{eqDo2} is equivalent to the so-called time-change formula derived  in \cite{BP} for multivariate regularly varying rf's.
}

Next, since for stationary $X$  we have that \eqref{convY} holds,  then in view of \cite{BP,klem} $X$ is  a multivariate regularly varying rf
 and $Y$ is the so-called  tail  rf of $X$, whereas   $\Theta$ is the so-called spectral tail rf. \EHF{Therefore for a stationary max-stable rf $X$ the  rf $\Theta$ defined in \eqref{eqtheta} is simply the spectral tail rf of $X$.}

Adopting the terminology of \cite{kab2009}  for stationary max-stable rf's $X$, 
we shall refer to  their spectral rf's $Z$ as  Brown-Resnick stationary (abbreviated  as BRs) rf's.\\
 \EHF{From $Z$ we can easily define the spectral tail  rf $\Theta$. Moreover, as mentioned in \eqref{fle} \cL{we simply} have $\Theta \equaldis Z$ if $Z(0)=1$  almost surely.  }
\EHF{The key properties of BRs   rf's $Z$ and spectral tail rf's $\Theta$ are the  TSF \eqref{eqDo} and the identity \eqref{eqDo2}, respectively. This is revealed by our next result}, which shows how to construct a BRs rf $Z$ from  a given rf $\Theta^*$ that satisfies  \eqref{eqDo2} and $\Theta^*(0)=1$ almost surely,   extending thus   \cite{Janssen}[Thm 4.2]  to rf's.
 
Let in the following $$\IA(p \cdot Y)= \min( i\in \TT :  \max_{j\in \TT}\abs{p_j Y(j)} = \abs{p_i Y(i)}),$$ where $p_j's$ are non-negative  numbers such that 
$\sum_{j\in \TT} p_j^\alpha =1$ (recall $\alpha=1$ in our case).\\
Hereafter $N$ is a rv independent of any other random element such that 
$\pk{ N=j}= p_j>0, j\in \TT$.  
 Further,  both $\min$ and $\max$ are defined with respect to a translation-invariant order on $\Z^d$, see \cite{BP} for the definition.

\BEL \label{prop2} 
 If $Y(t)= R \Theta^*(t),t\in \TT$ with $R$ a unit Pareto rv independent of $\Theta^*$ which satisfies \eqref{eqDo2} and $\Theta^*(0)=1$ almost surely, then $Z_N$ given by 
\bqn{ 
	Z_N(t) = \frac{ B^N Y (t)}{\max_{i\in \TT} \cL{p_i}  \EHF{B^N} Y(i) } \ind{ \IA( p \cdot B^N  Y ) =N}  , \quad t\in \TT
	\label{sm}
}
is a spectral rf \EHF{of some stationary max-stable  rf $X(t),t\in \TT$ with unit \FRE marginals}. \EHF{Moreover, the spectral tail rf $\Theta$ of $X$   has the same law as $\Theta^*$.}
\EEL

\BRM  \label{remF} \EHF{i) \cL{In case $\alpha\not=1$ the above construction is still valid by substituting the denominator with 
		$ (\max_{i\in \TT} {p_i}^\alpha  \EHF{B^N} Y(i)  )^{1/\alpha}  $.} }	
In fact,   \eqref{sm} is a  minor modification of the construction given in \cite{klem}[Prop 2.12]. The other known constructions in \cite{Janssen,Hrovje, klem} can be easily extended for the case $d>1$, we omit the details. \\
ii) A $\R^q$-valued  rf $\Theta(t),t\in \TT$ is called a spectral tail rf if it satisfies \eqref{eqDo2} where $\Theta(h),\Theta(-h)$ are substituted by $\norm{\Theta(h)},\norm{\Theta(-h)}$  with $\norm{\cdot}$ a norm on $\R^q$ and $F$ is redefined accordingly and further $\pk{  \norm{\Theta(0)}=1}=1$, see e.g., \cite{BojanS,Hrovje,BP}. For such a rf,  
a BRs rf $Z_N$ can be determined as in 
\eqref{sm} by changing  $\sum_{t\in \TT} p_t B^N Y(t)$ to 
$\sum_{t\in \TT} p_t B^N\norm{Y(t)}$ and instead of $\max_{t\in \TT} p_t B^N Y(t)$ and  $p \cdot B^N Y$  putting  
$\max_{t\in \TT} p_t B^N \norm{Y(t)}, p \cdot B^N \norm{Y(t)}$,  respectively (with $Y(t)= R \Theta(t)$ and $R$ a unit Pareto rv independent of $\Theta$). 
\ERM

\section{Classical, block \& functional  indices}
As mentioned in the Introduction the classical extremal index $\tX$ of a stationary max-stable  rf $X $  always exists.  
We show first that it is equal to the block extremal index $\widetilde{\tX}$ defined in \eqref{bloc} and then answer the question when $\tX=0$. This is already known for  $d=1$, see \cite{debicki2017approximation}. Our main findings  in \cL{\netheo{thK}} gives several formulas for $\tX$. 
\cL{The next result is a minor generalisation of the case $d=1$ stated in \cite{MR3745388}.}

\BEL \label{nub} If $X(t),t\in \TT$ is a stationary max-stable rf, then  $\tX = \cL{\widetilde{\tX}}$. 
\EEL

\noindent
Below we slightly modify  the definition of anchoring maps  introduced in \cite{BP}. 
Write next $\bar \Z^d$ for $ \Z^d  \cup\{ \IF\}$ and recall that $E=[0,\IF)^{\TT} $ is equipped with the product $\sigma$-field $\mathcal{E}$. 

\def\BTT{{\bar \Z}^d}
\def\So{O}
\BD We call a measurable map  $\mathcal{I}:  E \mapsto \BTT $ 
anchoring  if for \EHF{$\So=\{ f\in E: \mathcal{I}(f) \in \TT\}$}  
the following conditions are satisfied for  all $ f \in \So, i\in \TT$:\\
$i)$ $\mathcal{I}(f) =i$ implies $f(i) \ge \cL{\min( f(0) ,1) }$; \\
$ii)$ $ \mathcal{I}(f)=\mathcal{I}(B^i f) -i$. 
\ED

As in  \cite{BP} we define two important anchoring maps which \cL{are specified with respect to a translation-invariant  order on $\Z^d$. In 
	particular the  minimum and maximum below are with respect to such an order. An instance of a translation-invariant 
	order is  the lexicographical one.}  \EHF{  Hereafter  $\mathcal{S}(f)= \sum_{t\in \TT} f^\alpha(t)$ for any $f\in E$. Note that apart from Section 5.2  we have considered for simplicity only the case  $\alpha=1$.} \\

\noindent
{\bf Example 1.} Let the  non-empty set $\So  \in \mathcal{E}$ be given by 
$$ \So=\Bigl \{ f \in  E : \mathcal{S}( f)< \IF,\quad \max_{i\in \TT} f(i) > 0  \Bigr \}
$$
and define the first \EHF{maximum} functional 
$$
\IA(f)= \min \Bigr(j\in \TT :  f(j)= \max_{i\in \TT} f(i)\Bigl), \quad f\in \So,
$$
where  $\IA(f)=  \IF$ if $f\not \in \So$. 
Clearly, $\IA(f)$ is finite for $f\in \So$ and condition $i)$  holds by the definition, whereas  condition $ii)$ follows by the invariance (in the sense of \cite{WS}) of the translation-invariant order. 

The first and last maximum functionals are   important since they are both anchoring and $0$-homogeneous.  Moreover, for a stationary max-stable  rf $X(t),t\in \TT$ with spectral rf $\Theta$ and \FRE marginals $\Phi(x)= e^{-1/x^\alpha}, x>0$ we have that the law of $X$ is specified by  $\IA$ and $\Theta$ as follows 
\bqn{ - \ln  \pk{ X(i) \le x_i, i\in \TT } 
	&=& \sum_{i\in \TT} \frac 1 {x_i^\alpha } \pk{ \IA(\Theta/ (B^{-i} x)) =0}
	\label{nathana}
}	
for any $x= (x_i)_{i\in \TT}$ with finitely many positive components and the rest equal to \EHF{$\IF$}; here 
$\Theta/(B^{-i} x)= (\Theta(j)/ x_{j+i})_{j\in \TT}.$  The  proof of \eqref{nathana} is displayed   in Appendix,  see also \cite{Htilt}[Eq. (6.10)]. \cL{Note in passing that \eqref{nathana}  shows that the law of $X$ is uniquely determined by $\Theta$.} \\

\bigskip

\noindent
{\bf Example 2.} Define the first exceedance functional  by 
$$ \IAfE(f)= \min \Bigl(j\in \TT :  f(j)>1  \Bigr), \quad f \in  O$$
and  set 	 $\IAfE(f)= \IF$ if $f\not \in O$,  where 
$$\So=\Bigl \{ f \in  E: \mathcal{S}( f) < \IF, \quad {\max_{t\in \TT} f(t) }> 1 \Bigr \}  \in \mathcal{E}.  
$$
Clearly,  $\IAfE(f)$ for $f\in \So$ is finite and $i)$ holds. Moreover since $\IAfE(f), f\in \So$ is determined by a finite number of points  \EHF{in a neighbourhood of 0, then $\IAfE$ is  measurable.}  Again condition $ii)$ is implied by  the translation-invariance of the chosen 
 order on $\Z^d$. 

\cL{We call a measurable map $F: E \mapsto [0,\IF]$ shift-invariant if $F(B^h f)=F(f), h \in \TT, f\in E$.}

\BEL Let  $\Theta(t),t\in \R^d$ be a real-valued rf satisfying \eqref{eqDo2} with $\Theta(0)=1$ almost surely. If  $R$ is a unit Pareto rv  independent of $ \Theta$, then for any two anchoring maps  $\mathcal{I},\mathcal{I}'$  and \cL{any shift-invariant  map $F$} we have  (set  $ Y(t)= R  \Theta(t),t\in \TT$)
\bqn{ \pk{ \mathcal{I}(  Y  ) =0,  \mathcal{I}'(  Y ) \in \TT, \cL{F(Y)}<\IF} 
	= 
	\pk{ \mathcal{I}'(  Y) =0  , \mathcal{I}(  Y) \in \TT, \cL{F(Y)}<\IF}. \quad 
	\label{shum01}
}  
Moreover, $	\pk{\mathcal{I}(  Y) =0, \cL{F(Y)}< \IF }=0$ is equivalent with   $	\pk{\mathcal{I}(  Y) \in \TT, \cL{F(Y)} < \IF}=0.$

\label{ereLazri}
\EEL 

\BRM  If   $\mathcal{I}(Y),\mathcal{I}'(Y)$ are  almost surely in $\Z^d$,  then \eqref{shum01} boils down to  
$\pk{ \mathcal{I}'(  Y) =0}=\pk{ \mathcal{I}(  Y) =0}$, which is  already shown in \cite{BP}[Lem 3.5].   
\EHF{ In general, $\mathcal{I}(Y)$ might not be finite almost surely.}
\COM{ Therefore the event $\{\mathcal{S}(  Y) < \IF \}$ enters in our calculations.  In fact, under the conditions of \CE{\netheo{thK}} below on both $\mathcal{I}, \mathcal{I}'$ we shall show that 
	\bqn{ 
		\tX =\pk{ \mathcal{I}'(  Y) =0,\mathcal{S}(  Y) < \IF}=\pk{ \mathcal{I}(  Y) =0,\mathcal{S}(  Y) < \IF }.
	}	 
}
\ERM 

\noindent
Hereafter we consider    anchoring  maps  $\mathcal{I}:  E \mapsto \bar{\Z}^d$  such that  
\bqn{\label{vap}
	\pk{ \mathcal{I}(Y)\in \TT , \mathcal{S}(Y)< \IF}= \pk{\mathcal{S}(Y)< \IF},
}
which is in particular valid for both first (last) \EHF{maximum}  and first (last) exceedance functionals.  

\BEL \label{PropA}
If  $X(t),t\in \TT$ is  a stationary max-stable rf with some spectral rf $Z$ and spectral tail rf 
$\Theta$, then $\tX=0$  if and only if 
$\pk{\SST = \IF}=\pk{\SSZ = \IF}=1$.   If  further  the  anchoring map 
$\mathcal{I}$ 
satisfies \eqref{vap},  then $\tX=0$ is equivalent with 
\bqn{ \label{hofnung2} \pk{ \mathcal{I}(Y)= 0 , \mathcal{S}(Y)< \IF}=0.
} 
\EEL

Since the first and last maximum functionals are  0-homogeneous and finite on the set $\So=\{f \in E: \mathcal{S}(f)  < \IF ,\max_{i\in \TT} f(i) > 0\}$  we have  that $\pk{\SSZ = \IF} =1$ is equivalent with 
$$ \pk{\IA(Z) \not \in \TT } =1$$
and the same also holds for the last \EHF{maximum} functional.\\ 
In view of   \nelem{PropA}, \nelem{lemshift} and \cite{dom2016}  $\tX=0$ is equivalent with  $\pk{\SSZ \cL{=} \IF}=1$. Further  we have the following equivalent statements  (below $\norm{\cdot}$ is a norm on $\R^d$):\\ 
{\bf A1}: $  Z(t) \to 0 $ almost surely as $\norm{t} \to \IF$;\\
{\bf A2}:  $  \Theta(t) \to 0 $ almost surely as $\norm{t} \to \IF$;\\
{\bf A3}:  $  \mathcal{S}(Z) < \IF $ almost surely;\\
{\bf A4}:  $  \mathcal{S}(\Theta) < \IF $ almost surely.\\
The equivalence of {\bf A1} and {\bf A3} is shown in \cite{dom2016}, whereas the equivalence of {\bf A1} and {\bf A2} \cL{is a direct consequence of \nelem{lemshift} and similarly for the equivalence of {\bf A3} and {\bf A4}}. \cL{The equivalence {\bf A2} and {\bf A4} follows from \cite{Janssen} and \cite{WS}. Note  further that 
	 $ Y(t) =R \Theta(t) \to 0 $ almost surely as $\norm{t} \to \IF$ is equivalent with {\bf A2} and 
	 $  \mathcal{S}(Y) = R   \mathcal{S}(\Theta) < \IF $ almost surely is equivalent with {\bf A4}. }

We state next the main result of this section;  \cL{define in the following 
$\mathcal{B}(Y)=  \sum_{ t\in \TT} \ind{Y(t)> 1}$} and 
interpret  $0:0$ and $\IF:\IF$ as  0.

\BT Let  $\mathcal{I},X$ be  as in 
\nelem{PropA}. If      $\mathcal{I} $ satisfies \eqref{vap} and $ \pk{\mathcal{S}(\Theta) < \IF}>0$,  then  
\bqn{ \tX &=& \pk{\mathcal{I}(Y) =0,\mathcal{S}(Y) < \IF} \label{lucsherohuA}\\
	&=& \pk{\IAfE(Y) =0 }  \label{BI2}  \\
	&=& \pk{\IA(\Theta) =0}  \label{BI1}  \\
	&=& \pk{\mathcal{I}(\Theta) =0,\mathcal{S}(\Theta) < \IF} \label{lucsherohu}\\
	& =& \E*{ \frac{\max_{t\in \TT} \Theta(t) }{ \sum_{t\in \TT} \Theta(t)} } \label{heng}\\
	& =& \E*{ Z(0)\frac{\max_{t\in \TT} Z(t) }{ \sum_{t\in \TT} Z(t)} }\\
	&=& 	  \E*{\frac{1}   {\mathcal{B}(Y)} },
	\label{lucejeluce}
\label{VH}
}
where \eqref{lucsherohu}  holds if further $\mathcal{I}$ is 0-homogeneous.  \cL{Moreover $ \{\mathcal{B}(Y) < \IF\} = 
\{\mathcal{S}(Y)< \IF\}$ almost surely} and in particular $\tX=1$ if and only if $\Theta(i)=0$ \EHF{almost surely} for all $ i\in \TT, i\not=0$.
\label{thK}
\ET 

\BRM \label{nms} i) $\Theta(t)= \Theta_1(t_1)\Theta_2(t_2), t_1 \in \Z^k, t_2 \in \Z^m, t=(t_1,t_2)\in \TT$ with $\Theta_1,\Theta_2$ independent rf's satisfying \eqref{eqDo2} and $\pk{ {\Theta_i(0)}=1}=1,i=1,2$, then  \eqref{heng} implies that  $\tX= \theta_{X_1} \theta_{X_2}$ where $X, X_i,i=1,2$ are stationary max-stable  rf's with spectral rf $\Theta$ and $\Theta_i,i=1,2$, respectively.  \\  
ii) For $d=1$ and $\tX=1$ the claim  that $\Theta(i)=0, i\not=0$ in \netheo{thK}  follows also from \cite{MR3454028}[Prop 2.2 (ii)].\\
iii) \EHF{Since $\Theta$ uniquely defines $X$,  then \netheo{thK} implies  that the only stationary max-stable  rf $X$ such that $\tX=1$ is that given by \eqref{lfever}. In view of \eqref{convT}  $\Theta(i)=0,i\not=0$ and hence by \eqref{sm} 
	$$Z_N(t) = \frac{ 1 }{ \cL{p_t}}\ind{N=t} , \quad t\in \TT$$
	is a spectral rf for $X$ specified in \eqref{lfever}, where $N$ is a discrete rv with positive probability mass function $p_t>0,t\in \TT$.}\\
iv) Taking $F(f)= \ind{\mathcal{I}(f)=0, \mathcal{S}(f) < \IF}$, then \eqref{lucsherohu} implies $\tX=\nu_{X,F}$ under the further assumption that $\mathcal{I}$ is a $0$-homogeneous functional satisfying \eqref{vap}. \\
v) \cL{It follows from the proof of \netheo{thK} that \eqref{lucejeluce} holds without the assumption that $\pk{\mathcal{S}(\Theta)< \IF}>0$. Hence $\theta_X=0$ if and only if $\mathcal{B}(Y) = \IF$ almost surely. Further, from \netheo{thK} we have that {\bf A1, A2, A3} and {\bf A4} are  equivalent  with 
{\bf A5}:  $ \mathcal{B}(Y) <\IF $ almost surely.} \\
iv) \cL{ Formula \eqref{heng} appears initially  as extremal index in \cite{Genna04,Genna04c} and  in \cite{DiekerY} as Pickands constant.}

\ERM

\section{The anti-clustering condition}
Since stationary max-stable rf's with \FRE marginals are multivariate regularly varying  (see for more details \cite{BP}) 
the classical extremal index  of those rf's can be calculated using the findings of \cite{BP} and \cite{WS}.  
In the framework of stationary multivariate regularly varying   rf's    
the  anti-clustering condition of \cite{Davis} plays a crucial role for the calculation of   extremal index. Considering  the stationary max-stable rf $X(t),t\in \TT$ with unit \FRE marginals, in view of \cite{BP} the aforementioned condition reads as follows:\\

\noindent
{\bf Condition C}: Suppose that there exists a positive sequence of non-decreasing integers $r_n \to \IF$ as $n\to \IF$ and $\limit{n} \cL{r_n^d/n}=0$ such that for any $s>0$ 
\bqny{ \limit{m} \limsup_{n \to \IF} \pk*{ \max_{ m < \norm{t} < r_n, t\in \Z^d  }  X(t)> ns \lvert X(0)> ns }=0.}

The equivalence of Condition C and $\pk{\SST < \IF}=1$ for the case $d=1$ is known, see \cite{klem}. The case $d\ge 1$ of  Brown-Resnick model  is dealt with in \cite{WS}[Prop 6.2]. Next we show that this equivalence  holds for a general stationary max-stable rf $X$ with spectral tail  rf $\Theta$ and spectral rf $Z$.  

\BEL  The anti-clustering Condition C for $X$ is equivalent with {\bf Ai}, $\mathbf{i}=1, ...,  5$.
\label{blok2}
\EEL    
If $\pk{\SST < \IF}=1$ or equivalently Condition C holds, then by \cite{BP} 
 \cL{\nelem{nub}}, \nelem{blok2}  and \cite{BP}\cL{[Prop 5.2]} 
for any anchoring map $\mathcal{I}$ 
\bqn{ \tX=\pk{\mathcal{I}(Y)=0}= \pk{\IA(Y)=0}= \pk{\IA(\Theta)=0}\in (0,1],
	\label{fito}
}
provided that $\pk{\mathcal{I}(Y) \in \TT}=1$. In the special case  $\mathcal{I}=\mathcal{I}_{fe}$ (as shown already in \cite{BP})  
\bqn{ \label{erdhi}
	\tX=\pk{ \max_{ 0 \prec t} Y(t) \le 1}.
}	
Here $\prec$ denotes a translation-invariant order on $\TT$.

\BRM \cL{ The expression in \eqref{erdhi} is a well-known formula in the Gaussian setup and has appeared in numerous papers inspired by \cite{Albin1990}. This special formula for the Gaussian setup is also referred to as Albin's constant, see \cite{DiekerY}.  In the context of  stationary regularly varying time series the same formula has appeared   in \cite{BojanS}.}
	\ERM 

\EHF{Next, consider the case that Condition C does not hold, i.e., 	$p=\pk{\SST< \IF}\in (0,1)$ and define the rf's ${\Theta}_1= \Theta  \lvert (\SST < \IF) $   and $\Theta_2=\Theta \lvert ( \SST = \IF)$.
	In view of  \cite{dom2016}[Thm 9, Prop 10], for two independent 
	stationary max-stable rf's $\eta_i(t), t\in \TT, i=1,2$ with unit \FRE marginals and corresponding  spectral tail rf's   equal in law  to $\Theta_i,i=1,2$ we have that   $X$ has the same law as 
\bqn{\label{sL}
	 \max( p \eta_1(t), (1- p) \eta_2(t)) , \quad t\in \TT.
} 
	Since   $\eta_1$ satisfies  Condition C, then by \cite{WS}[\cL{Prop 5.2}], \cL{\nelem{nub}},  \eqref{fito} and \netheo{thK}   
	\bqn{ \tX = p \pk{ \IA(\Theta_1)=0} = p \theta_{\eta_1} \in (0,1].
		\label{magix}
	}
	Alternatively, since  by the stationarity  of $X$ we have that $\tX$ exists and moreover  $\theta_{\eta_2}=0$,  then  \nelem{lemTogo} implies that $ \tX = p \theta_{\eta_1}$.  Consequently, we conclude that Condition C, \nelem{lemTogo}, \cL{representation \eqref{sL} }together with the findings of \cite{BP} establish  \cL{the validity of the  first four expressions}  in \netheo{thK}.\\
We remark that from the above arguments, by \eqref{erdhi} and  \nelem{theo1} we obtain  
\bqny{ 
	\tX		&=& \E*{ \max_{ 0 \preceq t }  \Theta(t)  -\max_{ 0 \prec t }  \Theta(t) ; 
		\cL{S(\Theta) < \IF}}\\
	& = &\E*{ \max_{ 0 \preceq t }  Z(t)  -\max_{ 0 \prec t }  Z(t); \cL{ \mathcal{S}(Z)< \IF}    }. 
} 	
The first formula above is already obtained for   the  Brown-Resnick model (see Section \ref{sekA}) in \cite{WS}[\cL{Corr 6.3}] and for the case $d=1$ in \cite{MR3745388}[Thm 2.1]. 
}

\section{Examples} \label{sekA}
We present below some  examples starting first with the   Brown-Resnick model. The second example  and  \nelem{prop2}  show in particular   how to construct stationary max-stable rf's starting from any  $\alpha$-summable deterministic   sequence. \EHF{We then discuss how to construct from some given rf  a stationary max-stable rf $X$ such that $\tX$ equals a   given constant.}
 
\subsection{Brown-Resnick model}
Consider $Z(t)= e^{ W(t)- \sigma^2(t)/2}, t\in \TT$ with $W(t),t\in \TT$ a centered Gaussian rf with   variance function $\sigma^2$ which is not identical to 0 and $\sigma(0)=0$. Let $X(t),t\in \TT$ denote a max-stable rf with spectral rf $Z$. The case 
 $W$ is a standard Brownian motion and $d=1$ is investigated in \cite{bro1977} and therefore this construction is referred to as the Brown-Resnick model.\\
For any fixed $h\in \TT$ the Gaussian rf (set $ \gamma(s,t)= Var(W(t)- W(s)),s,t\in \TT$)
$$S_h(t)= W(t)- W(h)-\gamma(h,t)/2, \quad  \forall t\in \TT
$$    is such that $S_h(h)=0$ almost surely  and has variance function
$\sigma_h^2(t)= \gamma(h,t)$.\\
With the same arguments as in \cite{Htilt}, it follows that   $Z_h(t)=e^{S_h(t)},t\in \TT$ is also a spectral rf for  $X$ for any $h\in \TT$. Since $S_h(t), t\in\Z^d$ is a \cL{Gaussian} rf with variance $Var(W(t)- W(h))= \gamma(t,h)$, then  the law of $X$ depends only on $\gamma(h,t)$ and not on $\sigma^2$. 
If we assume that  $W$ has   stationary  increments, then  \eqref{eqDo} implies  that  $X$ is  a stationary max-stable  rf. 
The fact that  $Z_h(h)=1$ for any $h\in \TT$ almost surely implies that $\Theta:=\Theta_0$ defined in \eqref{eqtheta} is simply given by   $\Theta(t)= Z(t),t\in \TT$ and hence  (recall $Y= R \Theta$) 
$$Y(t)=e^{\widetilde W(t) + Q}, \quad \widetilde W(t)= W(t)- \sigma^2(t)/2, \quad t\in \TT, 
 $$
where $Q= \ln R$ is a unit exponential rv independent of $W$.\\
 For an $N(0,1)$ rv $V$ with distribution $\Phi$  being  independent of $Q$  and all $c>0,x\inr $ 
(set  $\bar \Phi=1- \Phi, V_c=cV  - c^2/2$)
\bqn{   \pk{ V_c+ Q > x }
	 &=&   \cL{ \pk{ V_c + Q > x,   V_c> x }+
		\pk{ V_c + Q > x, V_c\le x  }}\notag\\
	 &=&   \cL{ \pk{ V_c > x  }+ e^{-x} \E{ e^{ V_c} \ind{ V_c\le x }}} \notag\\
	&=&  
		   \pk{ V_c > x } + e^{-x} \pk{ cV \le x- c^2/2},  
	\label{gg1}
}
where we used that the exponentially tilted rv $U$ defined by $\pk{ U \le x}= 
\E{ e^{ V_c} \ind{ V_c\le x }}, x\inr$ has $N(c^2/2,c^2)$ distribution, see e.g., \cite{Htilt}[Lem 7.1]. Consequently, for all  $ t\in \TT$ such that $c:=\sigma(t)>0$ and all  $ y>0$
\bqn{ \quad  \quad 
	\pk{ Y(t) \le y} 
	&=& \Phi( c^{-1} \ln y+ c/2) - e^{-\cL{1/y}} \Phi( c^{-1}  \ln y- c/2), \quad 
}	
which  agrees with the claim of  \cite{WS}[Prop 6.1] where  the stationary case  is considered.

Next, under the assumption that $W$ has stationary increments,  in view of  \eqref{heng} and \eqref{lucejeluce} 
\COM{\bqny{ \pk{ cW + E - c^2/2 > x } &=& \int_{\R} \pk{ E> x+c^2- cw} \varphi(w) dw\\
		&=& \int_{cw < x+c^2/2}  e^{w - c^2/2- x}  \varphi(w) dw + \pk{cW > x+c^ 2/2}\\
		&=& e^{- c^2/2 - x + 1/2}(2 \pi)^{-1/2} \int_{cw < x+c^2/2}  e^{ -(  w-1)^2/2 }   dw + \pk{W > x+a}\\
		&=& e^{- a - x + 1/2}(2 \pi)^{-1/2} \int_{w < x+a-1}  e^{ -w^2/2 }   dw + \pk{W > x+a}\\
		&=& e^{- a - x + 1/2} \pk{ W< x+a - 1}  + \pk{W > x+a}\\
	}
	\bqn{ \pk{ cW + E - c^2/2 > x } &=& \int_{0}^\IF  \pk{ -cW + c^2/2+x<  t}e^{- t} dt \\
		&=& \int_{cw < x+c^2/2}  e^{w - c^2/2- x}  \varphi(w) dw + \pk{cW > x+c^ 2/2}\\
	}
}
\bqn{ \label{aldous} 
	\tX 
	&=& \E*{ \frac{1}{ \sum_{t\in \TT} \ind{ \widetilde W(t) +Q >0}}   }=
	\E*{ \frac{\max_{t\in \TT} e^{ \widetilde W(t)} }{ \sum_{t\in \TT} e^{ \widetilde W(t)} }  } ,
	\label{skb}
} 
 \cL{which yields the following lower bound}
\bqn{
	\tX 	&=& \E*{ \frac{1}{ \sum_{t\in \TT} \ind{ \widetilde W(t) +Q >0}}   } 
	 \ge \frac{1}{ \E{  \sum_{t\in \TT} \ind{\widetilde  W(t)  +Q>0}}}\notag \\
	&=&\frac{1}{ \sum_{t\in \TT} \pk{ \widetilde  W(t) +Q >0} }\notag \\
	&=&\frac{1}{ \sum_{t\in \TT} \bar \Phi( \sigma^2(t)/2)}, 
	\label{xho}
}	
where we used Fubini theorem for the first equality and \eqref{gg1} implies  \eqref{xho}. The lower bound above is strictly positive under some growth conditions on $\sigma$, see \cite{MR4029237} for similar calculations in the continuous case. \cL{Derivation of a tight positive lower bound is of general interest since in most of the cases direct evaluation of $\theta_X$ is not feasible.}  \\

It is of some interest to compare two different extremal indices of stationary max-stable Brown-Resnick rf's for different variance functions. With similar arguments as in \cite{debicki2017approximation}[Thm 3.1] we can prove the following result:

\BEL Let $X_1(t), t\in \TT$ and $X_2(t), t\in \TT$ be two  stationary max-stable Brown-Resnick  rf's corresponding to two  centered Gaussian processes $W_1, W_2$ with stationary increments, continuous trajectories  and variance functions $\sigma_1^2$ and $\sigma_2^2$ which vanish at the origin.  If 
$\sigma_1(t)\ge \sigma_2(t)$ holds for all $t\in \TT$,   then $\theta_{X_1} \ge \theta_{X_2}$.
\label{LL}
\EEL

\BRM 
i) Under the conditions of \nelem{LL} 
$$ \E*{ \frac{1}{ \sum_{t\in \TT} \ind{ \widetilde W_1(t) +Q >0}} }
\ge \E*{ \frac{1}{ \sum_{t\in \TT} \ind{ \widetilde W_2(t) +Q >0}} }.
$$
ii) The calculation of $\tX$ and different expressions for it have appeared in the literature in various contexts: the most prominent one concerns  extremes of Gaussian rf's  where in fact $\widetilde{\tX}$ has been originally calculated, see e.g., \cite{PickandsB,debicki2002ruin, Pit19}.  
The first expression in  \cL{\eqref{aldous} for the continuous setup, $d=1$ and the fractional Brownian motion case 
is obtained  in \cite{Berman92}[Thm 10.5.1].}   Applications to sequential analysis and statistics  have given rise to various forms of formula \eqref{skb}, see e.g., \cite{Sigm1,Sigm2}.  As already shown in \cite{DiekerY}  \eqref{skb}  is useful for simulations of $\tX$. 
\ERM
\subsection{$\Theta$ generated by summable sequences} 
Let $c_i,i\in \TT$ be non-negative constants  satisfying  $\sum_{i\in \TT} c_i^\alpha =C \in (0,\IF) $ for some $\alpha >0$ and define 
$$ \Theta(i)= \frac{c_{i+ S}}{c_S}, \quad i\in \TT$$
for a given rv $S$ with values in $\TT$ satisfying  
$$\pk{S=i}= c_i^\alpha/ C, \quad i\in \TT .$$
Clearly, $\Theta(0)=1$ almost surely and moreover  $\Theta$ satisfies  \eqref{eqDo2} stated for the case $\alpha>0$ as below, namely  for any $h\in \TT$
\bqny{
	\E*{{  \Theta^\alpha (h)} F({ \Theta})} &=& \E{ c_{h+S}^\alpha/c_{S}^\alpha \ind{ c_S \not=0}F(c_{\cdot + S}) }\\   
	&=& \frac{1}{ C}\sum_{i\in \TT}   c_{h+i}^\alpha  \ind{c_i\not=0}   F(c_{\cdot + i}) \\
	& =&  \E*{  F( B^{h} {\Theta }) \ind{ \Theta(-h)\not=0 } }
}
is  valid for any $0$-homogeneous  measurable functional $F:  E \mapsto [0,\IF]$.\\ 
Clearly, $\SST=\sum_{t \in \TT} \Theta^\alpha(t)$ is finite almost surely, hence 
\bqn{\tX = \E*{ \frac{ \max_{t\in \TT} c_{t+ S}^\alpha   }{ \sum_{t\in \TT} c_{t+ S}^\alpha   }}
	= \frac{ 1} C  \max_{t\in \TT} c_{t}^\alpha  \in (0,1].
	\label{need}
}
\EHF{We note that $\tX$ given in \eqref{need} is the extremal index of a large class of stationary rf's, see e.g., \cite{MR3226001,SojaA}. 
}
\subsection{Constructions of $X$ with  given extremal index} 
\EHF{From the previous example we conclude that for any $a \in (0,1]$ we can construct a stationary max-stable  rf $X$ such that $\tX=a$. We present below  examples of  rf $X$ satisfying  $\tX=0$ and then we construct stationary max-stable rf's $X^{(p)}$ indexed by $p\in (0,1)$  and calculate their extremal indices.\\  
	Consider next independent,  non-negative rf's $\Theta_k(t), t\in \Z, k\le d $ that satisfy  \eqref{eqDo2} such that 
	$\pk{\Theta_k(0)=1}=1, k\le d$. It follows that the rf  $\Theta(t) = \prod_{ 1 \le k \le d} \Theta_k(t_k), t= (t_1 \ldot t_k) \in \TT$ also satisfies \eqref{eqDo2}.  In view of \nelem{prop2} we can construct  stationary max-stable rf's  $X, X_k, k\le d$  corresponding to $\Theta, \Theta_k, k\le d$.  As already mentioned in Remark \ref{nms} $ii)$ we have $\tX= \prod_{k\le d} \theta_{X_k}$ and therefore $\tX=0$ if some $\theta_{X_k}$ equals zero. If we  define 
	$\Theta_k(j)=1$ for all even integers $j$ and $\Theta_k(j)=0$ for all odd integers $j$, then   $\Theta_k$ satisfies \eqref{eqDo2}.  Since  $\mathcal{S}(\Theta_k) = \IF$ almost surely, then $\theta_{X_k}=0$ follows and hence also $\tX=0$. \\   
	In view of our examples, we can construct two independent stationary max-stable rf's $\eta_1(t), \eta_2(t), t\in \TT$ with unit $\FRE$ marginals and spectral tail rf's $Z_1$ and $Z_2$, respectively satisfying $\pk{\mathcal{S}(Z_1)< \IF}= \pk{\mathcal{S}(Z_2)= \IF}=1$. 
	The rf $X^{(p)}(t)= \max(p \eta_1(t), (1- p) \eta_2(t)), t\in \TT$ for any given $p\in (0,1)$ is stationary and further max-stable with unit \FRE marginals. As already shown in the previous section, we have $\theta_{X^{(p)}}= p \theta_{\eta_1}$.
}
\section{Proofs}  

\prooflem{theo1}
For a given non-negative spectral rf $Z$ of a max-stable rf $X$ with unit \FRE marginals  by the de Haan representation of $X$ 
for any  $t_i \in \Z^d, x_i \in (0,\IF), i\le n$
\bqn{\label{fidis} -\ln \pk{ X(t_1) \le x_1 \ldot X(t_n) \le x_n}= \E*{ \max_{1 \le i \le n} \frac{Z(t_i)}{x_i}}.}
Consequently, with  \EHF{ $t_0=h\in \TT $ and $x_0=1$} we obtain as $u\to \IF$ 
\bqny{
	\lefteqn{
		\pk{ u^{-1} X(t_i)\le x_i, i=1 \ldot n \lvert X(t_0)> u} }\\
	&\sim & 
	u \pk{ u^{-1} X(t_i)\le x_i, i=1 \ldot n, u^{-1}X(t_0)> x_0}\\
	& = &	u [\pk{ u^{-1} X(t_i)\le x_i, i=1 \ldot n} - \pk{ u^{-1} X(t_i)\le x_i, i=0 \ldot n}]\\
	&\to &  \E*{    \max_{i=0 \ldot n} \frac{Z(t_i)}{x_i}-   \max_{i=1 \ldot n} \frac{Z(t_i)}{x_i}}, \quad u\to \IF \\\
	&= &  \E*{   \mathbb{I}(Z(t_0)>0)\Bigl  [ \max_{i=0 \ldot n}\frac{ Z(t_i)}{x_i} -    \max_{i=1 \ldot n} \frac{Z(t_i)}{x_i}\Bigr] }\\
	&= &  \E*{   Z(t_0) \mathbb{I}(Z(t_0)>0) \Bigl [ \max_{i=0 \ldot n} \frac{Z(t_i)}{Z(t_0)x_i} -    \max_{i=1 \ldot n}\frac{ Z(t_i)}{Z(t_0) x_i} \Bigr] }\\
	&= &  \E*{    \max_{i=0 \ldot n}  \frac{ \Theta_h(t_i)}{x_i} -    \max_{i=1 \ldot n} \frac{\Theta_h(t_i)}{x_i} },
}  
where the last line follows by the definition of $\Theta_h$ in \eqref{eqtheta}.  \EHF{Hence in view of \eqref{convY} and the fact that 
	$\Theta_h(h)=1$ almost surely,}  the proof is complete.
\QED

\bigskip 
\noindent 
\prooflem{prop2}  
{Since by the assumptions $\sum_{j\in \TT}p_j =1$  and $\Theta^*$ is non-negative we have  \EHF{for any $j\in \TT$}
$$ \E*{ \sum_{i\in \TT} \EHF{p_i}\Theta^*(i-j)}= \sum_{i\in \TT} p_i \E{ \Theta^*(i-j)}=\sum_{i\in \TT} p_i \pk{ \Theta^*(j-i) >0}\le 1,$$
which together with the non-negativity of $\Theta^*$ implies for some norm $\norm{\cdot}$ on $\R^d$
\bqn{\label{reta} \lim_{ \norm{t} \to \IF, t\in \TT} p_t \Theta^*(t\EHF{-j})=\lim_{ \norm{t} \to \IF, t\in \TT} p_t Y(t\EHF{-j}) =0}
almost surely. Consequently, since further  
$$ \pk{ p_N> 0}=\pk{Y(0)>1}=1,
$$ 
then    $\max_{t\in \TT} p_t \EHF{B^N} Y(t) \in (0, \IF) $ almost surely and thus  $Z_N$ in \eqref{sm} is well-defined. 
}
Next, for any $\EHF{a}, h\in \TT$ and any $0$-homogeneous measurable functional $F:  E \mapsto [0,\IF]$, by the independence of $N$ and $Y$ applying Fubini theorem  we obtain  
\bqny{ 
	\lefteqn{ \E{ Z_N(h) F( \EHF{B^a} Z_N)}}\\
	 &=& \E*{ \frac{ B^N Y (h)}{\max_{s\in \TT} p_s \EHF{B^N} Y(s) } \ind{ \IA( p \cdot B^N  Y ) =N} F(B^{a+N} Y)} \\
	&=& \sum_{j\in \TT} \E*{ p_j\frac{ B^j \Theta^* (h)}{\max_{s\in \TT} p_s   \Theta^*(s\EHF{-j}) } 
		\ind{ \IA( p \cdot B^{j}  \Theta^* ) =j} F(B^{a+j} \Theta^*)}\\
	&=& \sum_{j\in \TT} \E{ B^j \Theta^* (h)  \ind{ \IA( p \cdot B^j  \Theta^* ) =j} F(B^{a+j} \Theta^*)}  \\ 
	&=& \sum_{j\in \TT} \E{   \ind{ \IA( p \cdot B^h  \Theta^* ) =j, \Theta^*(j-h) >0} F(B^{a+h} \Theta^*)}  \\ 
	&=& \E*{F(B^{a+h} \Theta^*)  \sum_{j\in \TT}  \ind{ \IA( p \cdot B^h  \Theta^* ) =j, \Theta^*(j-h) >0} }  \\ 
	&=& \E{F(B^{a+h} \Theta^*) }\\
	&=& \E{Z_N(a) F(B^h Z_N) },	
}
where the third equality follows since  	$\IA( p \cdot B^j  \Theta^* ) =j$ implies 
$$ \EHF{\max_{ s\in \TT} p_s \Theta^*(s-j)}=  p_j B^j \Theta^*(j) = p_j\Theta^*(0)=\EHF{p_j}>0$$
almost surely, the fourth equality follows from  \eqref{eqDo2} and \EHF{the assumption that $\pk{\Theta^*(0)=1}=1$}, the sixth one is consequence of the following (which follows from \eqref{reta})
\bqny{\label{peta}
	\sum_{j\in \TT}  \ind{ \IA( p \cdot B^h  \Theta^* ) =j}=\ind{ \IA( p \cdot B^h  \Theta^* ) \in \TT}=1
}	
almost surely and the fact that $\IA( p \cdot B^h  \Theta^* ) =j$ implies \EHF{for any $h\in \TT$}
$$ \EHF{p_j} \Theta^*(\EHF{j}-h) \ge p_h \Theta^*(0) \ge p_h>0$$
almost surely \EHF{and consequently  $\Theta^*(\EHF{j}-h)>0$ almost surely}. 
Finally,   the last claimed equality is established  by repeating the calculations for $\E{Z_N(a)F(B^h Z_N)}$. Hence the proof follows   by \eqref{eqDo} and the definition  of the  spectral tail rf $\Theta$ via the spectral rf $Z$. \QED

\bigskip 
\noindent 
\prooflem{nub} Let $r_n\in \TT,n\ge1$ be non-negative integers with components $r_{nj},j\le d$ such that $\limit{n} n/r_{nj}=\limit{n} r_{nj}=\IF$. 
The stationarity of $X$ yields further  
$$ C(A)= \E*{ \max_{ i\in A } Z(i)} = C(A')$$
for any finite  set of indices $A \subset \TT$ 
and any $A'\subset \TT$ which is a shift/translation  of $A$. Moreover,  by the sub-additivity of the maximum
$$C(A\cup B) \le C(A)+ C(B).
$$
Hence the growth of $C(A)$ is as that of the counting measure of $A$, see \cite{DM} for this argument and \cite{MR534842}. Consequently, 
\bqny{ 
	\limit{n} \frac{\E*{ \max_{ 0 \le   i \le r_n, i\in \TT } Z(i)}}{ \prod_{j=1}^d r_{nj}}  
	&=& \limit{n} n^{-d}{\E*{ \max_{ i \in [0,n]^d , i\in \TT } Z(i)}} = \mathcal{H}.
}
The assumption on $r_n$ and \eqref{fidis} imply that 
\bqny{ 
	\widetilde{\tX} \sim \frac{\pk{ \max_{ 0 \le   i \le r_n, i\in \TT } X(i)> n } }{ \prod_{j=1}^d r_{nj}  
		\pk{ X(0)> n }}	\sim  \frac{\E*{ \max_{ 0 \le   i \le r_n, i\in \TT } Z(i)}}{ \prod_{j=1}^d r_{nj}}, \quad n\to \IF.
}
Hence  $\mathcal{H}= \tX$ establishes the proof. \QED 

\bigskip 
\noindent 
\prooflem{ereLazri} We give first a   key characterisation of tail rf's  proved initially in \cite{Hrovje} and also 
stated for rf's in \cite{BP}. Namely, for any  measurable map $F:  E \mapsto [0,\IF]$
\bqn{ \label{tiltY}
	\E*{ F( Y) \ind{ {Y(i)}> 1/t}}&=&  t\E*{ F(B^{i}  Y) \ind{ { Y(-i)}> t}}
}	
holds for all $i\in \TT, t>0$.	If $\mathcal{I}, \mathcal{I}'$ are two  anchoring  maps, since $ {  Y(0)}= R> 1$ almost surely and $\mathcal{I}(Y)=i$ implies $Y(i)  >1$ almost surely,  by  \eqref{tiltY}
\bqny{ 
	\lefteqn{ \pk{\mathcal{I}(Y) \in \TT, \mathcal{I}'( Y)=0, \cL{F(Y)}< \IF } }\\
	&=& 	
	\sum_{i\in \TT} \pk{\mathcal{I}( Y )=i,\mathcal{I}'( Y)=0, \cL{F(Y)}< \IF }\notag \\
	&=& 	\sum_{i\in \TT} \pk{\mathcal{I}( Y)=i, { Y(i)}  > 1, \mathcal{I}'(  Y)=0 ,\cL{F(Y)}< \IF} \notag \\
	&=& 	\sum_{i\in \TT} \pk{\mathcal{I}(B^{i}  Y) =i, {  Y(-i)}> 1,\mathcal{I}'(B^{i}  Y) =0 ,\cL{F(Y)}< \IF} \notag \\
	&=& 	\sum_{i\in \TT} \pk{\mathcal{I}(  Y) =0, \cL{F(Y)}< \IF, { Y(-i)}> 1,\mathcal{I}'(  Y) =-i } \notag \\
	&=& 	\sum_{i\in \TT} \pk{\mathcal{I}(  Y) =0, \cL{F(Y)}< \IF, \mathcal{I}'(  Y) =-i } \notag \\
	&=& 	\pk{ \mathcal{I}'(  Y) \in \TT,  \mathcal{I}(  Y) =0 ,\cL{F(Y)}< \IF}.
}
With similar arguments we obtain  
\bqny{ \pk{ \mathcal{I}( Y)\in \TT,\cL{F(Y)}< \IF } 
	&=& 	\sum_{i\in \TT} \pk{\mathcal{I}( Y) =0,  \cL{F(Y)}< \IF,  { Y(-i)}> 1 }.
}	
Consequently, 
$	\pk{\mathcal{I}(  Y) =0,\cL{F(Y)}< \IF }=0$ 
is equivalent with   
$$	\pk{\mathcal{I}(  Y) \in \TT, \cL{F(Y)}< \IF }=0$$ 
establishing the proof. \QED

\bigskip 
\noindent 
\prooflem{PropA}   
As shown in \cite{dom2016} condition  $\pk{\SSZ = \IF }=1$ is equivalent with $X$ being  generated by a non-singular conservative flow. The latter is equivalent with $\tX=0$, see \cite{MR2453345} (which follows by  \cite{Genna04} if $d=1$  and by  \cite{MR2384479} for $d>1$). \cL{In view of \nelem{ereLazri} and \eqref{vap}
$ \pk{ \mathcal{I}(Y)=0 , \mathcal{S}(Y)< \IF} =0$ is equivalent with $\pk{\mathcal{S}(Y)< \IF} =0$. Applying 
\nelem{lemshift} in Appendix the latter is equivalent with $\pk{\SSZ<\IF}=0$. This establishes the proof since the latter is equivalent with  $\theta_X=0$.} \QED 

\bigskip 
\noindent 
\noindent 
\prooftheo{thK} We have that $\pk{\SSZ \EHF{<} \IF}=0$ is equivalent with $X$ is generated by a non-singular conservative flow, which in view of \cite{Genna04,MR2384479,Roy1} is equivalent with $\tX=0$. Applying  \nelem{ik} in Appendix  to BRs spectral rf $Z$ we have  that 
$Z F(Z)$ is also a BRs spectral rf for any measurable functional $F: E \mapsto  [0,\IF]$, which is $0$-homogeneous and shift-invariant. Since  both $\mathbb{I} ( \mathcal{S}(f)  = \IF) , \mathbb{I} (\mathcal{S}(f)   < \IF), f \in E$ are measurable 0-homogeneous and shift-invariant functionals and by the above 
$$ \limit{n} \frac 1 {n^d} \E*{  \max_{t\in [0,n]^d \cap \TT} Z(t)\mathbb{I} (\SSZ  = \IF)}=0
$$
we have using further \eqref{AVI}   
\bqn{\label{schlaf} 
	\tX=\HWD &=& 
	\limit{n} \frac 1 {n^d} \E*{  \max_{t\in [0,n]^d \cap \TT} Z(t) } \notag \\
	&=&	 \limit{n} \frac 1 {n^d} \E*{  \max_{t\in [0,n]^d \cap \TT} Z(t) \mathbb{I} (\SSZ  < \IF)}.
}
Next, assuming  that $\pk{\SSZ< \IF} >0$ by \nelem{lemshift}  $\pk{ \SST < \IF}>0$ and the converse also holds.  
\cL{Setting $Z_*(t)= Z(t)\mathbb{I} (\SSZ  < \IF)$ by \nelem{ik} it is BRs and further  $\mathcal{S}(Z_*)< \IF$ almost surely.  In view of \nelem{lemX} we can assume that $\mathcal{S}(Z_*)>0$ almost surely.   Applying   \eqref{eqDo}  and using 
	the equivalence  of {\bf A1} and {\bf  A3}} we obtain further
\bqny{
	\tX 
	&=&	 \limit{n} \frac 1 {n^d}  \sum_{ h \in  [0,n]^d \cap \TT} \E*{  Z_*(h)  \frac{ \max_{t\in [0,n]^d \cap \TT} Z_*(t)}
		{ \sum_{ t \in  [0,n]^d \cap \TT}  Z_*(t)}}\\
	&=&	 \limit{n} \frac 1 {n^d}  \sum_{ h \in  [0,n]^d \cap \TT} \E*{  Z_*(0)  \frac{ \max_{t\in [0,n]^d \cap \TT}  B^hZ_*(t)}
		{ \sum_{ t \in  [0,n]^d \cap \TT}  B^hZ_*(t)} }\\
	&=&
	\cL{	 \lim_{\varepsilon \downarrow 0}  \limit{n} \frac 1 {n^d}  \sum_{ h \in  [\varepsilon n , (1- \varepsilon )n]^d \cap \TT} \E*{  Z_*(0)  \frac{ \max_{t\in [0,n]^d \cap \TT}  B^hZ_*(t)}
	{ \sum_{ t \in  [0,n]^d \cap \TT}  B^hZ_*(t)} }
}\\
	&=&	\E*{  Z_*(0)  \frac{ \max_{t \in \TT}  Z_*(t)}
		{ \sum_{ t \in  \TT}  Z_*(t)} } \\
	&=&	\E*{   \frac{ \max_{t\in \TT}  \Theta(t)}
		{ \SST }\mathbb{I} (\SST  < \IF)}.	
} 
Since by definition the events $\{ \IA(\Theta) \in \TT\}$ and 
$ \{ \SST < \IF\} $ are almost surely the same,    the $0$-homogeneity of $\IA(\cdot) $  implies  (recall 
$\Theta(0)=1$ almost surely)
\bqny{\tX
	&=& \E*{ \frac{\max_{t\in \TT} \Theta(t)}{\SST}\ind{\IA(\Theta) \in \TT  }}\\
	&=& \sum_{j \in \TT} \E*{ \frac{\max_{t\in \TT} \Theta(t)}{\SST} \ind{\IA(\Theta) = j   }}\\
	&=& \sum_{j \in \TT} \E*{  \Theta(j)\frac{\Theta(0)}{\SST} \ind{\IA(\Theta)  =j   }}\\
	&=& \sum_{j \in \TT} \E*{  \frac{\Theta(-j)}{\SST} \ind{\IA(B^j \Theta) =j   }}\\
	&=&  \E*{ \sum_{j \in \TT} \frac{\Theta(-j)}{\SST} \ind{ \IA(\Theta)   =0  } }\\
	&=& \pk{ \IA(\Theta)  =0}\\
	&=& \pk{ \IA(\Theta)  =0,\SST < \IF },
}
where we applied \eqref{eqDo2} in the  last third line combined with condition $ii)$ in the definition of anchoring maps  and also used that $\mathcal{S}(f),f\in E$ is a shift-invariant functional. 
Clearly, the last two formulas hold also for the last maximum functional. 
Since \eqref{vap}  implies 
\bqn{\label{hajv}
	\pk{ \mathcal{I}(Y) \not \in \TT ,\mathcal{S}(Y) < \IF}=0,
}
then using  \cL{\nelem{ereLazri}} to obtain the second equality below  
we have 
\bqny{
 	\pk{ \IA (\Theta)  =0, \SST < \IF} 
	&=&  \pk*{\IA (Y)  =0, \mathcal{S}(Y) < \IF, \mathcal{I}(Y) \in \TT }\\
	&& +  \pk*{\IA (Y)  =0, \mathcal{S}(Y) < \IF, \mathcal{I}(Y) \not \in \TT }\\
	&=&  \pk*{\IA (Y)  \in \TT , \mathcal{S}(Y) < \IF, \mathcal{I}(Y) =0  }\\
	&=&  \pk*{ \mathcal{I}(Y) =0 ,\mathcal{S}(Y) < \IF }
}
and hence $\tX= \pk{ \IAfE(Y) =0}$ follows and the same is true also for the last exeedance functional.
\cL{In view of the equivalence {\bf A2} and {\bf A4} we have 
\bqn{ \label{masha} 
	\{ \mathcal{S}(Y) < \IF\} \subset   
 \{ \mathcal{B}(Y)  < \IF\},
}
  with $ \mathcal{B}(Y):=\sum_{ t\in \TT} \ind{Y(t)> 1}$. Hence 
since $Y(0) =R \Theta(0) = R >1$ almost surely implies $\mathcal{B}(Y) \ge  1$ almost surely}
\bqny{    
\lefteqn{	\E*{\frac{ \mathcal{B}(Y) } { \mathcal{B}(Y) } \ind{ \JJY =0, \mathcal{S}(Y)< \IF }} }\\
	&=&  
\sum_{ t\in \TT}\E*{\frac{1}   {\mathcal{B}(Y) } \ind{ \JJY =0, Y(t)> 1, \mathcal{S}(Y)< \IF }}\\
	&=&  \E*{\frac{1}   {\mathcal{B}(Y) } \sum_{ t\in \TT}\ind{ \JJY=-t, Y(-t)> 1, \mathcal{S}(Y)< \IF }}\\
	&=&  \E*{\frac{1}   { \mathcal{B}(Y) } \ind{ \JJY \in \TT, \mathcal{S}(Y)< \IF }}\\
	&=&  \E*{\frac{1}   {\mathcal{B}(Y) } \ind{ \mathcal{S}(Y)< \IF }}\\
	&=&  \E*{\frac{1}   { \mathcal{B}(Y)} \ind{  \mathcal{S}(\Theta)< \IF }},
}
where we used \eqref{tiltY} to derive  the last fourth line and  the last second equality follows from \eqref{hajv}. 
\cL{With the same arguments as in the proof of \cite{PH2020}[Lem 2.5] considering the discrete setup as in \cite{klem}  for any $n>0$
$$\E*{  \max_{t\in [0,n]^d \cap \TT} Z(t) }   = 
 \sum_{t\in [0,n]^d \cap \TT}  \E*{ \frac{ 1 }{  \sum_{s\in [0,n]^d \cap \TT}  
 		 \ind{{Y(s-t)}> 1 } }}. $$
Since  $Y(0)>1$ almost surely and thus the denominator in the expectation above  is greater equal 1 and converges as $n\to \IF$ almost surely to 
$\mathcal{B}(Y)$, it follows by the dominated convergence theorem  that 
\bqny{  
	\tX &=&  \limit{n} n^{-d} \E*{  \max_{t\in [0,n]^d \cap \TT} Z(t) } = \E*{\frac{1}   { \mathcal{B}(Y)}  }\le 1,
}
hence \eqref{lucejeluce} holds. From the last two expressions of  $\tX$ we conclude that  
$\E*{\frac{1}   { \mathcal{B}(Y)} \ind{  \mathcal{S}(Y)= \IF }}=0$. Consequently,  almost surely 
$\{ \mathcal{B}(Y) <\IF \} \subset \{ \mathcal{S}(Y)< \IF \}$,   
which together with  \eqref{masha} implies that almost surely 
$$	\{ \mathcal{B}(Y) < \IF\}  = \{ \mathcal{S}(Y) < \IF\}.$$
} 
Next, if $ \pk{\Theta(i)=0}=1$  for all $i\not=0, i\in \TT$, then  
$$ \tX = \E*{ \frac{\max_{t\in \TT} \Theta(t)}{\sum_{t\in \TT} \Theta(t)}\ind{\SST < \IF } } =1.$$
Conversely, if $\tX=1$, then necessarily $\pk{ \SST < \IF}=1$ and thus 
 $$ \tX=1= \E*{ \frac{\max_{t\in \TT} \Theta(t)}{\sum_{t\in \TT} \Theta(t)} }$$
implying that $\max_{t\in \TT} \Theta(t) =\sum_{t\in \TT} \Theta(t)$ almost surely. Taking $\mathcal{I}(f)=\IA(f)$ we have that $\tX=\pk{\JJT=0}=1$ implies that $\max_{t\in \TT} \Theta(t)=\Theta(0)=1$ almost surely and therefore   
$$\sum_{t\in \TT} \Theta(t) = 1+ \sum_{t\in \TT, t\not=0} \Theta(t) =1$$
almost surely. Consequently, (recall $\Theta(i)$'s are non-negative) $\pk{\Theta(i) =0}=1$  for all  $i\not=0, i\in \TT$
establishing the proof.
\QED  

\noindent 
\prooflem{blok2} For any $s>0$ and any non-decreasing sequence of integers $r_n,n\in \mathbb{N}$ tending to infinity such that $\limit{n}\cL{r_n^d}/n=0$ we have for any positive integer $m$ (recall $\E*{Z(t)}=1$ for any $t\in \TT$)
$$n^{-1}{\E*{\max_{ m < \norm{t} < r_n, t\in \Z^d  } Z(t)}}\le n^{-1} {\sum_{ m < \norm{t} < r_n, t\in \Z^d  } \E*{Z(t)}} \to 0, \quad n\to \IF,
$$
hence by \eqref{fidis} and the dominated convergence theorem  
\bqny{
	\lefteqn{ 
		1- \limit{n}  \pk{ \max_{ m < \norm{t} < r_n, t\in \Z^d  }  X(t)> ns \lvert X(0)> ns }
	}\\
	&=& 
	s\limit{n} n \pk{ \max_{ m < \norm{t} < r_n, t\in \Z^d  }  X(t)\le  ns,  X(0)> ns }\\
	&=&  \E*{ \max_{ m < \norm{t} < \IF , t\in \Z^d , t=0 }  Z(t)- \max_{ m < \norm{t} < \IF , t\in \Z^d  }  Z(t)}\\
	&=&  \E*{ \ind{Z(0)>0} \Bigl[ \max_{ m < \norm{t} < \IF , t\in \Z^d , t=0 }  Z(t)- \max_{ m < \norm{t} < \IF , t\in \Z^d  }  Z(t) 
		\Bigr]}\\ 
	&=&  \E*{ Z(0) \ind{Z(0)>0}\Bigl[ \max_{ m < \norm{t} < \IF , t\in \Z^d , t=0 }  \frac{Z(t)}{Z(0)}- \max_{ m < \norm{t} < \IF , t\in \Z^d  }  \frac{Z(t)}{Z(0)} 
		\Bigr]}\\
	&=&  \E*{ \Bigl(1-  \max_{ m < \norm{t} < \IF , t\in \Z^d  }  \Theta(t) \Bigr)_+} 
}	
for any positive integer $m$ (recall $\Theta(0)=1$ almost surely).  If  {\bf A1} holds, then by the dominated convergence theorem 
$$\limit{m} \E*{ \max_{ m < \norm{t} < \IF , t\in \Z^d , t=0 }  Z(t)- \max_{ m < \norm{t} < \IF , t\in \Z^d  }  Z(t)} = \E*{Z(0)} =1,$$ hence Condition C is satisfied. \\ 
Conversely, if Condition C is satisfied  for some sequence $r_n,n\ge 1$ of non-negative increasing integers, then 
by the above calculations  
\bqny{ \lefteqn{ 1- \limit{m}\limit{n} \pk{ \max_{ m < \norm{t} < r_n, t\in \Z^d  }  X(t)> ns \lvert X(0)> ns }}\\
	& =& 
	\limit{m} \E*{(1-  \max_{ m < \norm{t} < \IF , t\in \Z^d  }  \Theta(t))_+} =1
}
and thus almost surely as $m\to \IF$
$$\max_{ m < \norm{t} < \IF , t\in \Z^d  }  \Theta(t) \to 0.
$$
Consequently, by \nelem{limTT} \cL{in Appendix} condition {\bf A2} holds, hence the proof follows from Remark \ref{nms}.
\QED 

\section{Appendix}

For notational simplicity we consider the case $\alpha=1$ in the following. The results for $\alpha>0$ can be formulated with obvious modifications.
\BEL \label{lemX} 
If  $X(t), t\in \TT$ is  a max-stable rf with de Haan representation \eqref{eq1} and some spectral rf   $Z$ satisfying  $\E{Z(t)} \in (0,\IF)$ for all $t\in \TT$, then we can find a spectral rf   $Z_*$ for $X$ such that $\max_{t\in \TT}Z_*(t)> 0$ almost surely. 
\EEL 
\noindent 
\prooflem{lemX} Let $w_i, i\in \TT$ be positive constants such that
 $$ \E*{\sum_{i\in \TT} w_i  Z(i) } \in (0, \IF).$$  
 $w_i$'s exist since $\E{Z(i)} \in (0,\IF)$ for any $i\in \TT$.  
By the choice of $w_i$'s we have that 
$$M= \max_{i\in \TT} w_i Z (i)$$
is a non-negative rv and $a=\E{M }  \in (0,\IF)$. Let $Z_*(t),t\in \TT$ be a rf defined  by 
$$ \pk{ Z_* \in A} = \E{ M   \ind{ a Z/M \in A } /a}$$
for any measurable set $A \subset E$. Since by the above definition 
$$   \pk{ \max_{i\in \TT } w_i Z_*(i) =0} = \E{ M  \ind{ \max_{i\in \TT} w_i Z(i)/M =0 } /a}
= 0$$
it follows  that $\pk{\max_{i\in \TT }  Z_*(i) =0}=0$.
Moreover, for any $x_i \in (0,\IF), t_i \in \TT, i\le n $  
\bqny{ \lefteqn{ - \ln \pk{ X(t_1) \le x_1 \ldot X(t_n) \le x_n} }\\
	& =& \E{ \max_{1 \le i\le n} Z(t_i)/x_i   }\\
	& =& \E{ \ind{\max_{1 \le i\le n} Z(t_i)>0} \max_{1 \le i\le n} Z(t_i)/x_i }\\
	& =& \E{ M/a \ind{M>0}\ind{ \max_{1 \le i\le n} Z(t_i)>0} \max_{1 \le i\le n} aZ(t_i)/(M x_i) }\\	
	& =& \E{\ind{ \max_{1 \le i\le n} Z_*(t_i)>0} \max_{1 \le i\le n} Z_*(t_i)/ x_i }\\	
	& =& \E{\max_{1 \le i\le n} Z_*(t_i)/ x_i },	
}
where the third equality is a simple consequence of $\max_{1 \le i\le n} Z(t_i)>0$ implies $M>0$. Hence $Z_*$ is a spectral rf   for $X$. 
\cL{The calculations above show that we can define alternatively $Z_*(t)= \pk{\max_{s\in \TT} Z(s)> 0 } Z(t)
$ conditioned on $\max_{ s\in \TT} Z(s)> 0$, which was   suggested by the reviewer}. 
\QED

\bigskip 
\noindent 
Proof of \eqref{nathana}: As in the proof of \nelem{lemX}, we can assume without loss of generality that $Z$ is such that $\max_{t\in \TT} (Z(t)/x_t)>0$ almost surely  for any $x=(x_j)_{j\in \TT}$ a positive  sequence. Suppose for simplicity that $\alpha=1$ and let next  $x$ be a sequence with finite number of positive  elements and the rest equal to $\IF$ (we interpret  $a/\IF$ as $0$). Since further  $Z/x$ consists of zeros and finitely many positive numbers, then $\IA(Z/x) \in \TT$ almost surely.  
Consequently, by \eqref{fidis}, Fubini theorem and the fact that  $\IA(Z/x)  =j$ 
implies $ \max_{i\in \TT  }  (Z(t_i)/ x_i)= Z(j)/x_j$ almost surely  
\bqn{ - \ln \pk{ X(i) \le x_i,i\in \TT }  
	&=& \E{ \max_{i\in \TT }  Z(t_i) /x_i \ind{ \IA( Z /x)  \in \TT}}  \notag \\
	&=&  \sum_{j\in \TT} \E{  \max_{i\in \TT   }  Z(t_i)/ x_i \ind{ \IA(Z/x)  =j } }  \notag \\
	&=&  \sum_{j\in \TT } \frac 1 {x_j}\E{ Z(j) \IA(Z/x)  =j}  \notag \\
	&=&  \sum_{j\in \TT  } \frac 1 {x_j}\E{ Z(0) \IA( B^j Z/x)  =j}  \notag \\
	&=&  \sum_{j\in \TT } \frac 1 {x_j}\E{  Z(0)\ind{ \IA( (B^j Z/ x)/  Z(0))  =j}}  \notag \\
	&=&  \sum_{j\in \TT } \frac 1 {x_j}\pk{  \IA(B^{j}( \Theta /(B^{-j}x)) )  =j}  \notag \\
	&=&  \sum_{j\in \TT } \frac 1 {x_j}\pk{  \IA(   \Theta/(B^{-j}x))  =0}  \notag ,
}
where the fourth first equality follows from  \eqref{eqDo} and  the  last equality follows since $\IA$ is an anchoring map. \QED

\BEL \label{lemshift} 
Let $Z(t), t\in \TT$ be a BRs rf satisfying \eqref{elit}.   If $F:E  \mapsto [0,\IF]$ is a shift-invariant and $0$-homogeneous measurable map, then   
$\E*{ F(Z)} =0$ is equivalent with $\E*{F(\Theta)} =0$. If further $F$ is bounded by 1, then 
$  \E*{ F(Z)} =1$ is equivalent with $ \E*{ F(\Theta)} =1$.  
\EEL 

\noindent 
\prooflem{lemshift}  
By   the shift-invariance of $F$  and \eqref{eqDo} we have 
\bqny{ 0 &=& 
	\E*{ F( \Theta)  } =  \E*{  Z(0)F( Z/Z(0) ) }
	= \sum_{i\in \TT} \E*{  Z(0)F( B^{-i} Z)}\\
	&=&  \sum_{i\in \TT} \E*{  Z(i)F(Z)}
	\ge 	\E*{  \Bigl(\max_{i\in \TT}  Z(i) \Bigr) F(Z)},
}
hence since $Z$ is chosen such that  $ \max_{i\in \TT} Z(i)>0$ almost surely, then $ \E*{F(Z)}=0$ follows. If $\E*{F(Z)} =0$, 
then $F(Z)=0$ almost surely and thus 
$$	0= \E*{ Z(0) F(Z)}= \E*{ F(\Theta)}=0$$
 follows. 
Next,   $ \E*{ F(\Theta)} =1$ is the same as $\E*{ 1-  F(\Theta)} =0$, which is equivalent with $\E*{1- F(Z)} =0$ as shown above, establishing   the proof. \QED

\BEL \label{ik} If $F:E \mapsto [0,\IF]$ is a $0$-homogeneous measurable functional and $Z(t),t\in \TT$ is a BRs rf, then $Z_*=Z F(Z)$ is also a BRs rf, provided that $\E{Z_*(t_0)}\in(0,\IF)$ for some $t_0\in \TT$.
\EEL
\noindent 
\prooflem{ik} Using  \eqref{eqDo} we have that $\E{ Z_*(t)}= \E{Z_*(t_0)} \in (0,\IF)$ for any $t\in \TT$ and in particular  $\pk{F(Z)=0}<1$ and $\pk{F(Z)= \IF}=0$. Since $F$ is $0$-homogeneous, we have that $Z_*$ satisfies   \eqref{eqDo}, which is an equivalent condition for a spectral rf to be a BRs rf,  see \cite{Htilt}. 
\QED 

\BEL \label{limTT} If $V(t),t\in \TT$ is a non-negative  rf, then $\pk{\limit{\norm{t}}V(t) = 0}=1$  is equivalent with there exists a non-decreasing sequence of integers $r_n,n\ge 1$ that converge to infinite as $n\to \IF$  such that 
\bqn{\label{kela}
	\limit{m} \limsup_{n\to \IF} \pk{ \max_{ m \le \norm{t} \le r_n} V(t) > \delta}=0
}
is valid for any $\delta>0$. 
\EEL 
\noindent 
\prooflem{limTT} It is well-known that (see e.g.,  \cite{MR1458613}[A1.3]) 
$$\pk{\limit{\norm{t}}V(t) = 0}=1$$  if and only if 
for all large $m$ and any $\delta, \varepsilon $ positive
$$ 
\pk{\max_{ \norm{t} \ge m } V(t) > \delta}< \varepsilon, 
$$
which clearly implies \eqref{kela}. Assuming that the latter condition holds, then  for given $\delta, \varepsilon$ positive  there exists $N$ such that for all $m,n$ larger than $N$ we have $ \pk{ \max_{ m \le \norm{t} \le r_n} V(t) > \delta}< \varepsilon$. Since $\limit{n}r_n=\IF$, then  
$ \pk{ \max_{ m \le \norm{t} } V(t) > \delta}\le  \varepsilon$, hence the claim follows. 
\QED 

\EHF{\BEL Let $\eta_i(t),i=1,2, t\in \TT$ be two independent stationary rf's with unit \FRE marginal distributions. If  the  extremal indices of both $\eta_1$ and $\eta_2$ exist, then the  rf $X(t)= \max(p\eta_1(t), (1-p)\eta_2(t)), t\in \TT$ has for any $p\in (0,1)$   extremal index  $ \tX =  p \theta_{\eta_1} + (1-p)\ \theta_{\eta_2} \in [0,1].$
	\label{lemTogo}
	\EEL 	
	\noindent 
	\prooflem{lemTogo} By the independence of $\eta_1$ and $\eta_2$ we have that  $X$ is stationary with unit \FRE marginal distributions. In order to show   the claim  it suffices to prove  that $ \max_{ t\in [0, n]^d} X(t) /n^d $ converges in distribution as $n\to \IF$ to $ (p \theta_{\eta_1}+ (1-p) \theta_{\eta_2}) \xi ,$ where $\xi$ is a unit \FRE rv. As $n\to \IF$, by the assumptions  $\max_{ t\in [0, n]^d} \eta_i(t) /n^d$ converge for $i=1,2$ in distribution to $p_i  \theta_{\eta_i} \xi_i$ with $\xi_1,\xi_2$ two independent unit \FRE rv's and $p_1=1- p_2=p$. Since $\max( p_1  \theta_{\eta_1} \xi_1, 
	p_2  \theta_{\eta_2} \xi_2 )$ has the same df as $(p_1  \theta_{\eta_1}+ p_2  \theta_{\eta_2}) \xi$, the claim follows by the independence of $\eta_i$'s and  Slutsky's lemma.  \QED 
}
\section*{Acknowledgments} 
Partial 	support by SNSF Grants 200021-175752/1  and  200021-196888 is kindly acknowledged. 
\cL{I am in debt   to the reviewer and the Editor for several suggestions, comments and corrections which improved the manuscript significantly. }

\bibliographystyle{abbrv}

\bibliography{EEEA}
\end{document}